\documentclass[12 pt]{amsart}

\usepackage{amscd,amsmath,latexsym,amssymb,amscd}
\usepackage[mathscr]{euscript} 
\usepackage[all]{xy}
\usepackage{layout}
\usepackage{pb-diagram}
\usepackage{pstricks}
\usepackage{enumerate}
\usepackage{dsfont}

\newtheorem{theorem}{Theorem}[section]

\newtheorem{lemma}[theorem]{Lemma}
\newtheorem{proposition}[theorem]{Proposition}
\theoremstyle{definition}
\newtheorem{definition}[theorem]{Definition}

\theoremstyle{definition}

\newtheorem{property}[theorem]{Property}

\newpsobject{grilla}{psgrid}{subgriddiv=1,griddots=10,gridlabels=6pt}
\def\Proof{\medskip\noindent{\bf Proof: }}

\def\Z{\mathbb{Z}}

\def\Q{\mathbb{Q}}

\def\I{\mathbb{I}}

\def\O{\mathcal{O}}
\def\x{{\bf{x}}}
\def\y{{\bf{y}}}
\def\o{{\mathfrak{o}}}
\def\s{{\mathfrak{s}}}

\DeclareMathOperator{\Hom}{\textup{Hom}}
\DeclareMathOperator{\maj}{\textup{maj}}
\DeclareMathOperator{\fmaj}{\textup{fmaj}}

\begin{document}

\title[A basis for the diagonally signed-symmetric polynomials] 
{A basis for the diagonally signed-symmetric polynomials }

\author[J.~M.~G\'omez]{Jos\'e Manuel G\'omez}
\address{Department of Mathematics,
Johns Hopkins University, Baltimore, MD 21218, USA}
\email{jgomez@math.jhu.edu}

\begin{abstract}
Let $n\ge 1$ be an integer and let 
$B_{n}$ denote the hyperoctahedral group of rank $n$. 
The group $B_{n}$ acts on the polynomial ring 
$\Q[x_{1},\dots,x_{n},y_{1},\dots,y_{n}]$ by signed permutations 
simultaneously on both of the sets of variables $x_{1},\dots,x_{n}$ 
and $y_{1},\dots,y_{n}.$ The invariant ring
$M^{B_{n}}:=\Q[x_{1},\dots,x_{n},y_{1},\dots,y_{n}]^{B_{n}}$  
is the ring of diagonally signed-symmetric polynomials. 
In this article we provide an explicit free basis of $M^{B_{n}}$ as a 
module over the 
ring of symmetric polynomials on both of the sets of variables 
$x_{1}^{2},\dots, x^{2}_{n}$ and  $y_{1}^{2},\dots, y^{2}_{n}$ 
using signed descent monomials. 
\end{abstract}

\maketitle

\section{Introduction}  

Let $V$ be an $n$-dimensional vector space over  a field 
$k$ of characteristic zero.  Suppose  that $W$ is  a
finite reflection group in $V$; that is, $W$ is finite subgroup 
of $GL(V)$  generated by elements 
of order $2$ that fix a hyperplane pointwise. 
Then $W$ acts  by ring automorphisms on the symmetric 
algebra $S(V^{*})$, 
where $V^{*}$ is the dual of $V$. If we give $V$ a basis, then 
$S(V^{*})$ can be identified with a polynomial ring 
$k[\x]:=k[x_{1},\dots, x_{n}]$. 
The action of the group $W$ on the polynomial ring  
$k[\x]$, under the above identification, has been 
classically studied. For example, by \cite[Theorem A]{Chevalley} 
the ring $k[\x]^{W}$ consisting of all $W$-invariant polynomials  
is itself a  polynomial ring  on $n$ homogeneous generators.  
Consider now the diagonal action of $W$ on the  symmetric algebra 
$S(V^{*}\oplus V^{*})$. If we give $V$ a basis as before, 
then $S(V^{*}\oplus V^{*})$ can be identified with a polynomial 
algebra $k[\x,\y]:=k[x_{1},\dots, x_{n},y_{1},\dots, y_{n}]$ and 
$W$ acts diagonally on it.  In this case the ring 
$M^{W}:=k[\x,\y]^{W}$ consisting of all diagonally 
$W$-invariant polynomials  is no longer 
a polynomial algebra. Note that the ring 
$R^{W}:=k[\x]^{W}\otimes k[\y]^{W}$  of all 
polynomials that are $W$-invariant in both of the sets of 
variables $\x$ and $\y$ is 
naturally a subring of $M^{W}$. Therefore we can 
see $M^{W}$ as a module over $R^{W}$. 
It can be seen that in fact $M^{W}$ is a  
free module over $R^{W}$ of rank $|W|$. 
This relies on the fact that $M^{W}$ is a Cohen-Macaulay 
ring which is true by \cite[Proposition 13]{HE}.   This article 
is concerned with the determination of explicit free bases 
of $M^{W}$ as a module over $R^{W}$ 
for a particular class of groups using elementary methods. 
For simplicity we work with rational coefficients although 
all the constructions 
provided here work for any field of characteristic zero.  

In \cite{Allen} Allen provided an explicit basis for the case 
of the symmetric group. More precisely, suppose that 
$W=\Sigma_{n}$ acts on the polynomial algebra 
$\Q[\x]=\Q[x_{1},\dots, x_{n}]$ by permutations of the 
variables $x_{1},\dots, x_{n}$. In this case the invariant 
ring $\Q[\x]^{\Sigma_{n}}$ is the ring of symmetric polynomials 
which is a polynomial algebra on the elementary symmetric 
polynomials. Let $\Sigma_{n}$ act diagonally on 
$\Q[\x,\y]:=\Q[x_{1},\dots, x_{n},y_{1},\dots, y_{n}]$. 
Then $M^{\Sigma_{n}}=\Q[\x,\y]^{\Sigma_{n}}$ is the ring of 
diagonally symmetric or multisymmetric polynomials. 
Given  $\pi\in \Sigma_{n}$ 
define the diagonal descent monomial
\[
e_{\pi}:=\prod_{i\in Des(\pi^{-1})}(x_{1}\cdots x_{i})
\prod_{j\in Des(\pi)}(y_{\pi(1)}\cdots y_{\pi(j)})
=\prod_{i=1}^{n}x_{i}^{d_{i}(\pi^{-1})}
 y_{i}^{d_{\pi^{-1}(i)}(\pi)},
\]
where $Des(\pi)$ denotes the descent set of $\pi$ and
$d_{i}(\pi^{-1})$, $d_{\pi^{1}(i)}(\pi)$ are integers
(see Section \ref{section 2} for the definitions). Then by 
\cite[Theorem 1.3]{Allen} the collection 
$\{\rho_{\Sigma_{n}}(e_{\pi}) \}_{\pi\in \Sigma_{n}}$ 
forms a free basis of $M^{\Sigma_{n}}$  as a module over 
$R^{\Sigma_{n}}=\Q[\x]^{\Sigma_{n}}\otimes \Q[\y]^{\Sigma_{n}}$, 
where $\rho_{\Sigma_{n}}$ is the averaging operator  
defined below.

The goal of this article is to show that an analogous 
construction works for the hyperoctahedral 
group $B_{n}$ acting on the polynomial algebra 
$\Q[\x]=\Q[x_{1},\dots, x_{n}]$ by signed permutations. 
In this case, it is easy to see that invariant 
ring $\Q[\x]^{B_{n}}$ consists of all symmetric 
polynomials on the variables $x_{1}^{2},\dots, x_{n}^{2}$. 
Suppose that $B_{n}$ acts diagonally 
on the polynomial ring 
$\Q[\x,\y]=\Q[x_{1},\dots, x_{n},y_{1},\dots,y_{n}]$ by 
signed permutations. Then the invariant 
ring $M^{B_{n}}=\Q[\x,\y]^{B_{n}}$ is the ring of diagonally  
signed-symmetric polynomials. A free 
basis of it as a module over 
$R^{B_{n}}=\Q[\x]^{B_{n}}\otimes \Q[\y]^{B_{n}}$  can be 
constructed in the same spirit as in the case 
of permutations. Given $\sigma\in B_{n}$,  define the 
diagonal signed descent monomial  
\[
c_{\sigma}:=\left(\prod_{i=1}^{n}x_{i}^{f_{i}(\sigma^{-1})} \right)
\left(\prod_{i=1}^{n}y_{|\sigma(i)|}^{f_{i}(\sigma)},
 \right).
\]
See Section \ref{section 3}  for the definition of the 
numbers $f_{i}(\sigma)$. The goal of this article 
is the following theorem.

\begin{theorem}
Suppose that $n\ge 1$. Then the collection 
$\{\rho(c_{\sigma})\}_{\sigma\in B_{n}}$ forms a free basis of 
$M^{B_{n}}$ as a module over 
$R^{B_{n}}$, where $\rho$ is 
the averaging operator.
\end{theorem}

A similar basis to the one given in the previous theorem was constructed 
in \cite{BB}. Moreover, in there  a nice combinatorial interpretation 
of the basis monomials was provided 
in terms of certain diagrams of the square lattice. The author would 
like to thank F. Bergeron and R. Biagioli for pointing out their work 
to him. 

\section{The symmetric group}\label{section 2}

In this section we provide a brief review  
of an explicit basis for the coinvariant ring for groups
of type $A$ using descent monomials constructed  
by Garsia and Stanton in \cite{GS}.  A construction 
of a free basis for the ring of diagonally symmetric polynomials 
as a module over the symmetric polynomials constructed 
by Allen in \cite{Allen} is also reviewed.

\subsection{Major index}
For every integer $n\ge 1$, let $\Sigma_{n}$ denote the symmetric 
group of self bijections of the set $\{1,2,\dots, n\} $.  
We use the notation $\pi=[\pi_{1},\dots, \pi_{n}]$ for an element 
$\pi\in \Sigma_{n}$  with $\pi_{i}=\pi(i)$ for $1\le i\le n$.  
Given $\pi\in \Sigma_{n}$ define its descent to be the set
\[
Des(\pi):=\{1\le i\le n-1 ~|~ \pi(i)>\pi(i+1)\}
\]
and for $1\le i \le n $ let
\[
d_{i}(\pi):=|\{j\in Des(\pi) ~|~ j\ge i \}|.
\]
The numbers $d_{i}(\pi)$ clearly 
satisfy the  following properties:
\begin{enumerate}
\item $d_{1}(\pi)\ge d_{2}(\pi)\ge\cdots \ge d_{n-1}(\pi)
\ge d_{n}(\pi)=0$, and
\item if $i<j$ and $d_{i}(\pi)=d_{j}(\pi)$ then
$\pi(i)<\pi(i+1)<\cdots<\pi(j)$.
\end{enumerate}
The major index of $\pi\in \Sigma_{n}$, denoted by $\maj(\pi)$, 
is defined to be 
\[
\maj(\pi):=\sum_{i=1}^{n}d_{i}(\pi)=\sum_{i\in Des(\pi)}i.
\]
In \cite{MacMahon}, MacMahon showed that this statistic 
is equidistributed with respect to the length function; that is,  
the number of permutations of length $n$ 
with $k$ inversions is the same as the number of 
permutations of length $n$ with major index equal to $k$.  
Note that the numbers 
$d_{1}(\pi)\ge d_{2}(\pi)\ge\cdots \ge d_{n-1}(\pi)$ 
are defined exactly to provide a partition of the integer 
$\maj(\pi)$.
  
\subsection{Descent monomials}
Suppose that $\x=\{x_{1},\dots, x_{n}\}$ is a set of algebraically 
independent commuting variables. Consider the polynomial 
algebra $\Q[\x]:=\Q[x_{1},\dots, x_{n}]$ seen as a graded ring 
with ${\rm{deg}}(x_{i})=1$  for $1\le i\le n$. The group 
$\Sigma_{n}$ acts naturally on the polynomial algebra 
$\Q[\x] $ by permuting the variables $x_{1},\dots, x_{n}$.
It is well known that the ring of $\Sigma_{n}$-invariants, 
$\Q[\x]^{\Sigma_{n}}$,  is a polynomial algebra on the generators 
$e_{1}(x_{1},\dots, x_{n}),\dots, 
e_{n}(x_{1},\dots, x_{n})$, where $e_{k}(x_{1},\dots, x_{n})$ 
is the $k$-th elementary symmetric polynomial 
\[
e_{k}(x_{1},\dots, x_{n})=\sum_{1\le i_{1}<\cdots< i_{k}\le n}
x_{i_{1}}x_{i_{2}}\cdots x_{i_{k}}.
\]
Suppose that $\pi\in \Sigma_{n}$. Define the descent 
monomial associated to $\pi$ to be
\[
a_{\pi}:=\prod_{i\in Des(\pi)}x_{\pi(1)}\cdots x_{\pi(i)}=
\prod_{i=1}^{n}x_{\pi(i)}^{d_{i}(\pi)}
=\prod_{i=1}^{n}x_{i}^{d_{\pi^{-1}(i)}(\pi)}.
\]

\noindent{\bf{Example:}} Suppose that $\pi=[6,2,1,4,3,5]$. Then 
$a_{\pi}=x_{1}x^{2}_{2}x_{4}x^{3}_{6}$.

\medskip

Note that  ${\rm{deg}}(a_{\pi})=\maj(\pi)$. 
In \cite{GS} Garsia and Stanton used 
Stanley--Reisner rings to show that these monomials provide 
a basis for the coinvariant algebra of type $A$. More precisely, 
let $I^{A}_{n}$ denote the ideal in $\Q[\x]$ generated by the 
symmetric polynomials $e_{1}(x_{1},\dots, x_{n}),\dots, 
e_{n}(x_{1},\dots, x_{n})$. Then $\Q[\x]/I^{A}_{n}$ 
is the coinvariant algebra  of type $A$. Let $\bar{a}_{\pi}$ 
denote the image of $a_{\pi}$ in the coinvariant algebra 
under the natural map.  In \cite{GS} it was proved that  
the collection $\{ \bar{a}_{\pi}\}_{\pi\in \Sigma_{n}}$ forms a basis of 
$\Q[\x]/I_{n}^{A}$ as a $\Q$-vector space. Moreover,  
the collection $\{a_{\pi}\}_{\pi\in \Sigma_{n}}$ 
provides a free basis for $\Q[\x]$ as a module over the 
symmetric polynomials $\Q[\x]^{\Sigma_{n}}$. This result 
has an interesting geometric application. Consider the flag manifold 
$U(n)/T$, where $T\subset U(n)$ is a maximal torus. Then 
$H^{*}(U(n)/T;\Q)$ can be identified with  the invariant algebra 
$\Q[\x]/I_{n}^{A}$, but under this identification 
we need to graduate the variables $x_{1},\dots, x_{n}$ with 
degree $2$. Therefore the descent monomials provide an 
explicit basis for the cohomology of the flag manifold $U(n)/T.$

\subsection{Diagonal descent monomials}
Let $\y=\{y_{1},\dots, y_{n}\}$ be another set of algebraically 
independent commuting variables of degree 1 
and consider the polynomial algebra 
$\Q[\x,\y]$. The 
symmetric group $\Sigma_{n}$  acts diagonally on this 
polynomial ring and the ring of $\Sigma_{n}$-invariants, 
$M^{\Sigma_{n}}:=\Q[\x,\y]^{\Sigma_{n}}$, is known as the ring of 
diagonally symmetric or multisymmetric polynomials. 
Note that the ring of polynomials that are symmetric in both 
the variables $\x$ and $\y$,  
$R^{\Sigma_{n}}:=\Q[\x]^{\Sigma_{n}}\otimes \Q[\y]^{\Sigma_{n}}$, 
is a subring of $M^{\Sigma_{n}}$ and thus $M^{\Sigma_{n}}$ can be 
seen as a module over $R^{\Sigma_{n}}$.  In \cite{Allen} Allen 
constructed a free basis for the module $M^{\Sigma_{n}}$ using a 
variation of the descent monomials. Given $\pi \in \Sigma_{n}$ 
define the diagonal descent monomial to be
\[
e_{\pi}:=\left(\prod_{i=1}^{n}x_{i}^{d_{i}(\pi^{-1})} \right)
\left(\prod_{i=1}^{n}y_{\pi(i)}^{d_{i}(\pi)}
 \right)=\prod_{i=1}^{n}x_{i}^{d_{i}(\pi^{-1})}
 y_{i}^{d_{\pi^{-1}(i)}(\pi)}.
\] 

\noindent{\bf{Example:}} Suppose that $\pi=[4,6,1,2,5,3]$. Then 
$e_{\pi}=x_{1}^{2}x_{2}^{2}x_{3}^{2}x_{4}x_{5}
y_{1}y_{2}y^{2}_{4}y_{5}y_{6}^{2}.$

\medskip

Note that the total degree 
of $e_{\pi}$ is given by  
${\rm{deg}}(e_{\pi})=\maj(\pi)+\maj(\pi^{-1})$ for every 
$\pi\in \Sigma_{n}$.  Consider the averaging operator
\begin{align*}
\rho_{\Sigma_{n}}:\Q[\x,\y]
&\to  \Q[\x,\y]^{\Sigma_{n}}=M^{\Sigma_{n}}\\
f&\mapsto \frac{1}{n!}\sum_{\pi\in \Sigma_{n}}
\pi\cdot f.
\end{align*} 
Thus by definition $\rho_{\Sigma_{n}}(e_{\pi})$ is a diagonally 
symmetric 
polynomial. By \cite[Theorem 1.3]{Allen} the collection 
$\{\rho_{\Sigma_{n}}(e_{\pi}) \}_{\pi\in \Sigma_{n}}$ forms a 
free basis of $M^{\Sigma_{n}}$ 
as a module over $R^{\Sigma_{n}}$. 
It turns out that this result 
also has an interesting geometric application. 
Let $B(2,U(n))$ be the geometric 
realization of the simplicial space obtained by 
defining $B_{k}(2,U(n))=\Hom(\Z^{k},U(n))$, where 
$\Hom(\Z^{k},U(n))$ is  the space of ordered commuting 
$k$-tuples in $U(n)$. In \cite{AG} it is proved that the diagonal 
descent monomials can be used to obtain an explicit basis of 
$H^{*}(B(2,U(n));\Q)$ seen as a module over $H^{*}(BU(n);\Q)$, 
where $BU(n)$ is  the classifying space of $U(n)$.

\section{The hyperoctahedral group}\label{section 3} 

In this section we provide analogue constructions to the ones 
presented in the previous section, where the symmetric group 
is replaced by the group of signed permutations. 

\subsection{Flag major index}
Suppose that $n\ge 1$ is an integer. Denote by 
$\I_{n}$ the set of integers between $-n$ and $n$ not 
including $0$; that is, 
\[
\I_{n}:=\{-n,-n+1,\dots,-1,1, \dots, n-1,,n\}.
\]
Let $B_{n}$ denote the group of bijections $\sigma:\I_{n}\to \I_{n}$ 
such that $\sigma(-k)=-\sigma(k)$ for all $k\in \I_{n}$, with the 
composition of functions as the group operation. Thus the group 
$B_{n}$ is the group of signed permutations,  also known 
as the hyperoctahedral group of rank $n$.  It is easy to see that 
$B_{n}$ is isomorphic to the semidirect product 
$\Sigma_{n}\ltimes (\Z/2)^{n}$. We use the following notation 
for elements   $\sigma\in B_{n}$.  Let $\sigma_{i}=\sigma(i)$ for 
$1\le  i\le n$, then 
we write $\sigma=[\sigma_{1},\dots, \sigma_{n}]$. 
Note that the group $B_{n}$ is the Weyl group associated to 
Lie groups of type $B_{n}$ and $C_{n}$ and that the symmetric 
group $\Sigma_{n}$ is naturally a subgroup of $B_{n}$. 
As in the case of the symmetric group, 
given $\sigma\in B_{n}$ define its descent to be the set
\[
Des(\sigma):=\{1\le i\le n-1 ~|~ \sigma(i)>\sigma(i+1)\}
\]
and for $1\le i \le n $ let
\[
d_{i}(\sigma):=|\{j\in Des(\sigma) ~|~ j\ge i \}|.
\]
As before the numbers $d_{i}(\sigma)$ satisfy the  following 
important properties:
\begin{enumerate}
\item $d_{1}(\sigma)\ge d_{2}(\sigma)\ge\cdots \ge d_{n-1}(\sigma)
\ge d_{n}(\sigma)=0$, and
\item if $i<j$ and $d_{i}(\sigma)=d_{j}(\sigma)$ then
$\sigma(i)<\sigma(i+1)<\cdots<\sigma(j)$.
\end{enumerate}
On the other hand, define
\begin{equation*}
\varepsilon_{i}(\sigma):= \left\{
\begin{array}{ccc}
0& \text{ if } &\sigma(i)>0,\\
1&\text{ if } &\sigma(i)<0,
\end{array}
\right. 
\end{equation*}
and 
\[
f_{i}(\sigma):=2d_{i}(\sigma)+\varepsilon_{i}(\sigma).
\]
It is easy to see that the numbers $f_{i}(\sigma)$ also satisfy 
the properties: 
\begin{enumerate}
\item $f_{1}(\sigma)\ge f_{2}(\sigma)\ge\cdots \ge f_{n}(\sigma)$, and
\item if $i<j$ and $f_{i}(\sigma)=f_{j}(\sigma)$ implies 
$\sigma(i)<\sigma(i+1)<\cdots<\sigma(j)$ and all of these numbers 
have the same sign.
\end{enumerate}
The flag major index of $\sigma\in B_{n}$, denoted by $\fmaj(\sigma)$ 
is defined to be 
\[
\fmaj(\sigma):=\sum_{i=1}^{n}f_{i}(\sigma)=
2\maj(\sigma)+\text{neg}(\sigma),
\]
where $\maj(\sigma)=\sum_{i\in Des(\sigma)}i$ is the major index of 
$\sigma$ and $\text{neg}(\sigma)=|\{1\le i\le n ~|~ \sigma(i)<0\}|$.
This statistic for elements in $B_{n}$ was introduced in \cite{AR} and 
further studied in \cite{ABR}, \cite{ABR1} as a generalization of the 
major index for the hyperoctahedral group. This tool 
has successfully  been used to study representation theoretical 
properties of the group $B_{n}$  (see for example \cite{ABR1}). 
Note that the numbers 
$f_{1}(\sigma)\ge \cdots \ge f_{n}(\sigma)$ provide a partition of 
the flag major index $\fmaj(\sigma)$ in a similar way as in the 
case of the major index of a permutation in $\Sigma_{n}$. 
Moreover, if $\sigma\in \Sigma_{n}$ then 
$\fmaj(\sigma)=2\maj(\sigma)$ so the flag major index 
is indeed a natural generalization of the major index.

\subsection{Signed descent monomials}
Suppose that $\x=\{x_{1},\dots, x_{n}\}$ is a set of algebraically 
independent commuting variables. Consider the polynomial 
algebra $\Q[\x]$ seen as a graded ring 
with ${\rm{deg}}(x_{i})=1$  for $1\le i\le n$. The group 
$B_{n}$ acts naturally on the polynomial algebra 
$\Q[\x] $ by degree preserving ring homomorphisms in the 
following way. If  $\sigma\in B_{n}$ then 
\[
\sigma\cdot (x_{1}^{i_{1}}\cdots x_{n}^{i_{n}})
:= \left(\frac{\sigma(1)}{|\sigma(1)|} \right)^{i_{1}}\cdots 
\left(\frac{\sigma(n)}{|\sigma(n)|} \right)^{i_{n}} 
x_{|\sigma(1)|}^{i_{1}}\cdots x_{|\sigma(n)|}^{i_{n}} 
\]
Thus each $\sigma$ permutes the variables $x_{1},\dots, x_{n}$ 
with a suitable sign change. It is easy to see that any polynomial 
in the ring of $B_{n}$-invariants, $\Q[\x]^{B_{n}}$, must be a 
symmetric polynomial in the variables $x_{1}^{2},\dots, x_{n}^{2}$. 
It follows  that $\Q[\x]^{B_{n}}$ is a polynomial algebra on the symmetric 
polynomials 
$e_{1}(x_{1}^{2},\dots, x_{n}^{2}),\dots, 
e_{n}(x_{1}^{2},\dots, x_{n}^{2})$. 
Suppose that $\sigma\in B_{n}$. Define the signed descent 
monomial
\[
b_{\sigma}:=\prod_{i=1}^{n}x_{|\sigma(i)|}^{f_{i}(\sigma)}
=\prod_{i=1}^{n}x_{i}^{f_{|\sigma^{-1}(i)|}(\sigma)}.
\]

\noindent{\bf{Example:}} Suppose that $\sigma=[-6,2,-1,-4,3,5]$. 
Then $b_{\sigma}=x_{1}^{3}x_{2}^{4}x_{4}x_{6}^{5}$.

\medskip

Note that ${\rm{deg}}(b_{\sigma})=\fmaj(\sigma)$ for every 
$\sigma\in B_{n}$.  The signed descent monomials can be used to 
obtain a basis for the 
coinvariant algebra for groups of type $B, C$. More precisely, let 
$I_{n}^{B}$ denote the ideal in $\Q[\x]$ generated by the 
elements $e_{1}(x_{1}^{2},\dots, x_{n}^{2}),\dots, 
e_{n}(x_{1}^{2},\dots, x_{n}^{2})$. Then the quotient 
$\Q[\x]/I^{B}_{n}(\x)$ is the coinvariant algebra in this case.  
Let $\bar{b}_{\sigma}$ denote the image 
of $b_{\sigma}$ in the coinvatiant algebra under the natural map. 
By \cite[Corollary 5.3]{ABR1}  the collection 
$\{ \bar{b}_{\sigma}\}_{\sigma\in B_{n}}$ forms a basis of 
$\Q[\x]/I^{B}_{n}(\x)$ as a $\Q$-vector space.  We can 
also see $\Q[x]$ as a module over 
the invariant ring $\Q[\x]^{B_{n}}$. Since 
$\{ \bar{b}_{\sigma}\}_{\sigma\in B_{n}}$ forms a basis of 
$\Q[\x]/I^{B}_{n}(\x)$ as a $\Q$,  then using \cite[Theorem 1.2]{Allen} 
it can be seen that $\{b_{\sigma}\}_{\sigma\in B_{n}}$ forms a free 
basis of $\Q[\x]$ as a module over $\Q[\x]^{B_{n}}$. This result has 
a geometric application as in the case of the symmetric group, 
namely, the signed descent monomials provide an explicit basis 
for the rational cohomology of the 
flag manifold $G/T$, for a compact connected Lie group $G$ of 
type $B_{n},C_{n}$ and a maximal torus $T\subset G$.

\subsection{Diagonal signed descent monomials}
Consider now $\y=\{y_{1},\dots, y_{n}\}$ another set of algebraically 
independent commuting variables of degree 1 
and consider the polynomial algebra 
$\Q[\x,\y]:=\Q[\x]\otimes\Q[\y]=
\Q[x_{1},\dots, x_{n},y_{1},\dots,y_{n}]$. The 
group $B_{n}$ acts diagonally on this polynomial ring; that is,  
$B_{n}$ acts as singed permutations simultaneously 
on the variables 
$x_{1},\dots,x_{n}$ and $y_{1},\dots,y_{n}$. Define 
$M^{B_{n}}:=\Q[\x,\y]^{B_{n}}.$ 
In other words $M^{B_{n}}$ is the ring of diagonally signed-symmetric 
polynomials.  Note that the ring of polynomials that are 
signed-symmetric on both the variables $\x$ and $\y$,
$R^{B_{n}}:=\Q[\x]^{B_{n}}\otimes \Q[\y]^{B_{n}}$, is 
a subring of $M^{B_{n}}$ and thus we can see $M^{B_{n}}$ as a module 
over $R^{B_{n}}$. As it was pointed out before $M^{B_{n}}$ 
is a free module over $R^{B_{n}}$ and the goal 
of this article is to construct an explicit basis for $M^{B_{n}}$ as a 
module over $R^{B_{n}}$. For this we will consider the 
following monomials.

\begin{definition}
Suppose that $\sigma\in B_{n}$. The diagonal signed descent 
monomial associated to $\sigma$ is defined to be 
\[
c_{\sigma}:=\left(\prod_{i=1}^{n}x_{i}^{f_{i}(\sigma^{-1})} \right)
\left(\prod_{i=1}^{n}y_{|\sigma(i)|}^{f_{i}(\sigma)}
 \right)=\prod_{i=1}^{n}x_{i}^{f_{i}(\sigma^{-1})}
 y_{i}^{f_{|\sigma^{-1}(i)|}(\sigma)}.
\]
\end{definition}

Note that for any $\sigma\in B_{n}$ we have 
${\rm{deg}}(\sigma)=\fmaj(\sigma)+\fmaj(\sigma^{-1})$.

\medskip

\noindent{\bf{Example:}} Suppose that $\sigma=[2,-1,-4,3]$. Then 
$c_{\sigma}=x^{3}_{1}x^{2}_{2}x^{2}_{3}x_{4}
y_{1}^{3}y_{2}^{4}y_{4}$.

\subsection{Averaging polynomials}
Consider the averaging operator
\begin{align*}
\rho:\Q[\x,\y]
&\to  \Q[\x,\y]^{B_{n}}=M^{B_{n}}\\
f&\mapsto \frac{1}{|B_{n}|}\sum_{\sigma\in B_{n}} 
\sigma\cdot f.
\end{align*} 
The map $\rho$ is a ring homomorphism that is surjective. 
Moreover, as a $\Q$-vector space $M^{B_{n}}$ is generated by 
elements of the form $\rho(m(\x,\y))$, 
where $m(\x,\y)$ is a monomial in $\Q[\x,\y]$.   
We will use the following notation. Suppose that $p=(p_{1},\dots,p_{n})$ 
is a sequence of non-negative integers. Then 
we write $x^{p}$ to denote the monomial 
$x_{1}^{p_{1}}\cdots x_{n}^{p_{n}}$.  

\begin{lemma}
Suppose that $p=(p_{1},\dots, p_{n})$ and $q=(q_{1},\dots, q_{n})$ 
are sequences of non-negative integers and let $m(\x,\y)=x^{p}y^{q}$. 
If $p_{k}+q_{k}$ is odd for some $1\le k\le n$ then $\rho(m(\x,\y))=0$.
\end{lemma}
\Proof
Suppose that $1\le k\le n$ is such that  $p_{k}+q_{k}$ is odd. Define 
\[
B_{n}^{+}=\{\sigma\in B_{n}~|~ \sigma(k)>0 \} \text{ and } 
B_{n}^{-}=\{\sigma\in B_{n}~|~ \sigma(k)<0 \}.
\] 
Note that 
$B_{n}=B_{n}^{+}\sqcup B_{n}^{-}$. Moreover, there is a bijection 
$\tau:B_{n}^{+}\to B_{n}^{-}$ defined by 
\begin{equation*}
\tau(\sigma)(i):= \left\{
\begin{array}{ccc}
\sigma(i)&\text{ if } &i\ne k,\\
-\sigma(i)& \text{ if } &i=k.
\end{array}
\right. 
\end{equation*}
By definition 
\[
\rho(m(\x,\y))=\frac{1}{|B_{n}|}\sum_{\sigma\in B_{n}}
c_{\sigma}x_{|\sigma(1)|}^{p_{1}}\cdots x_{|\sigma(n)|}^{p_{n}} 
y_{|\sigma(1)|}^{q_{1}}\cdots y_{|\sigma(n)|}^{q_{n}}, 
\]
where 
\[
c_{\sigma}=\left(\frac{\sigma(1)}{|\sigma(1)|} \right)^{p_{1}+q_{1}}\cdots 
\left(\frac{\sigma(n)}{|\sigma(n)|} \right)^{p_{n}+q_{n}}.
\]
For any $\sigma \in B_{n}^{+}$  we have 
$c_{\tau(\sigma)}=-c_{\sigma}$ since $i_{k}+j_{k}$ is odd. 
Therefore 
\begin{align*}
\rho(m(\x,\y))&=\frac{1}{|B_{n}|}\sum_{\sigma\in B_{n}}
c_{\sigma}x_{|\sigma(1)|}^{p_{1}}\cdots x_{|\sigma(n)|}^{p_{n}} 
y_{|\sigma(1)|}^{q_{1}}\cdots y_{|\sigma(n)|}^{q_{n}}\\
&=\frac{1}{|B_{n}|}\sum_{\sigma\in B_{n}^{+}}
(c_{\sigma}+c_{\tau(\sigma)})
x_{|\sigma(1)|}^{p_{1}}\cdots x_{|\sigma(n)|}^{p_{n}} 
y_{|\sigma(1)|}^{q_{1}}\cdots y_{|\sigma(n)|}^{q_{n}}=0.
\end{align*}
\qed

\medskip

By the previous lemma, it follows that $M^{B_{n}}$ is generated 
as a vector space over $\Q$  by the elements of the form 
$\rho(m(\x,\y))$, where 
$m(\x,\y)=x^{p}y^{q}$  and
$p=(p_{1},\dots, p_{n})$ and $q=(q_{1},\dots, q_{n})$ 
are  sequences of integers 
such that  $p_{k}+q_{k}$ is even for all 
$1\le k\le n$. Note that 
for any such monomial we have 
\[
\rho(m(\x,y))=\frac{1}{n!}\sum_{\alpha\in \Sigma_{n}}
x_{\alpha(1)}^{p_{1}}\cdots x_{\alpha(n)}^{p_{n}}
y_{\alpha(1)}^{q_{1}}\cdots y_{\alpha(n)}^{q_{n}}
\]
and thus $\rho(m(\x,\y))\ne 0$ for such monomials. 
Suppose that $\sigma\in B_{n}$ and consider the  
diagonal signed descent  monomial 
$c_{\sigma}=\prod_{i=1}^{n}x_{i}^{f_{i}(\sigma^{-1})}
 y_{i}^{f_{|\sigma^{-1}(i)|}(\sigma)}$ as defined above.
Using the properties of the numbers $f_{i}(\sigma)$, 
it is easy to see that for every $1\le i\le n$ the numbers 
$f_{i}(\sigma^{-1})$ and $f_{|\sigma^{-1}(i)|}(\sigma)$ have the 
same parity. Therefore by the previous comment we have 
$\rho(c_{\sigma})\ne 0$ for all $\sigma\in B_{n}$. By definition 
$\rho(c_{\sigma})\in M^{B_{n}}$. We will show below that 
the collection $\{ \rho(c_{\sigma})\}_{\sigma\in B_{n}}$ forms 
a free basis of $M^{B_{n}}$ as a module over $R^{B_{n}}$.

\subsection{Ordering of monomials} We will work with monomials 
$m(\x,\y)$ whose exponents are ordered in the following way.

\begin{definition} 
Suppose that $p=(p_{1},\dots, p_{n})$ and $q=(q_{1},\dots, q_{n})$ 
are two sequences of non-negative integers with 
$p_{k}+q_{k}$ even for all $1\le k\le n$. We say that the monomial
 $m(\x,\y)=x^{p}y^{q}$ is ordered and write $m(\x,\y)\in \O_{n}$ 
 if the exponents 
of $m(\x,\y)$ satisfy the following conditions:
\begin{enumerate}
\item $p_{1}\ge p_{2}\ge \cdots \ge p_{n}$,
\item if $p_{i}=p_{i+1}$ and they are even, then $q_{i}\ge q_{i+1}$, and  
\item if $p_{i}=p_{i+1}$ and they are odd, then $q_{i}\le q_{i+1}$.
\end{enumerate}
\end{definition}

The previous ordering can be described in the following way. For each 
integer $q$ define 
\begin{equation*}
\s(q):= \left\{
\begin{array}{ccc}
q&\text{ if } & q \text{ is even},\\
-q& \text{ if } &q \text{ is odd}.
\end{array}
\right. 
\end{equation*}
Then $m(\x,\y)=x^{p}y^{q}\in \O_{n}$ if and only if 
$(p_{1},\s(q_{1}))\ge_{\ell}\dots\ge_{\ell}(p_{n},\s(q_{n}))$, 
where $\ge_{\ell}$ denotes the lexicographic order. Using the 
definition given above it is easy to see that for every 
$\sigma\in B_{n}$ the diagonal signed descent monomial
 $c_{\sigma}$  is ordered in this way; that is, 
 $c_{\sigma}\in \O_{n}$ for all $\sigma\in B_{n}$.

Suppose now that $m(\x,\y)=x^{p}y^{q}$ is a monomial with 
$p_{k}+q_{k}$ even for all $1\le k\le n$ but whose exponents are 
not necessarily ordered as above. Consider the signed-symmetric 
polynomial $\rho(m(\x,\y))$. Note that in this polynomial there 
exists a unique monomial $n(\x,\y) $ whose exponents are 
ordered as above; that is,  $\rho(m(\x,\y))$ contains a unique monomial 
$n(\x,\y)\in \O_{n}$ and $\rho(m(\x,\y))=\rho(n(\x,\y))$. Because of this 
we can restrict ourselves to monomials $m(\x,\y)\in \O_{n}$. 

Next we define a total order on $\O_{n}$. 
For this suppose that $q=(q_{1},\dots, q_{n})$ is a sequence 
of integers. Define the ordering of $q$ to be
\[
\o(q):=(q_{\alpha(1)},\dots,q_{\alpha(n)}),
\]
where $(q_{\alpha(1)},\dots,q_{\alpha(n)})$ is a rearrangement 
of the sequence $q$  
in a decreasing way; that is,  $\alpha\in \Sigma_{n}$ is such that 
$q_{\alpha(1)}\ge\dots\ge q_{\alpha(n)}.$
Suppose now that  $m(\x,\y)=x^{p}y^{q}$ is a monomial 
with $p=(p_{1},\dots, p_{n})$ and $q=(q_{1},\dots, q_{n})$.
Define 
\[
\o(m(\x,\y)):=(\o(p),\o(q)).
\] 
In other words, $\o(m(\x,\y))$ 
recovers the exponents of the monomial $m(\x,\y)$ ordered 
in a decreasing fashion. Using this we can define 
the following total order on $\O_{n}$.

\begin{definition}
Suppose that $m(\x,\y)=x^{p}y^{q}$ and  
$n(\x,\y)=x^{i}y^{j}$ are  two monomials in $\O_{n}$.
We write $m(\x,\y)\succcurlyeq n(\x,\y)$ if and only if 
\begin{enumerate}
\item $\o(m(\x,\y))\ge_{\ell} \o(n(\x,\y))$, and 
\item if $\o(m(\x,\y))= \o(n(\x,\y))$ then 
$(p,\s(q))\ge_{\ell} (i,\s(j))$.
\end{enumerate}
\end{definition}
Here for a sequence of integers $q=(q_{1},\dots, q_{n})$ we 
have  $\s(q):=(\s(q_{1}),\dots, \s(q_{n}))$.  

\medskip

\noindent{\bf{Example:}} Suppose that $n=4$. Consider 
the monomials 
$m(\x,\y)=x^{7}_{1}x^{6}_{2}x^{6}_{3}x^{5}_{4}
y^{3}_{1}y^{8}_{2}y^{6}_{3}y^{5}_{4}$ and 
$n(\x,\y)=x^{7}_{1}x^{6}_{2}x^{6}_{3}x^{5}_{4}
y^{5}_{1}y^{8}_{2}y^{6}_{3}y^{3}_{4}$. Then 
$m(\x,\y), n(\x,\y)\in \O_{4}$  are such that
\[
\o(m(\x,\y))=\o(n(\x,\y))=(7,6,6,5,8,6,5,3)
\]
and $m(\x,y)\succcurlyeq n(\x,\y)$. 

\subsection{Signed index permutation} Suppose that 
$m(\x,\y)=x^{p}y^{q}\in\O_{n}$. Note that by 
construction the sequence $p=(p_{1},\dots, p_{n})$ is ordered 
in a decreasing way. This is not necessarily true for the sequence 
$q=(q_{1},\dots,q_{n})$. With this in mind, we can  associate to 
the monomial $m(\x,\y)$ the unique element 
$\sigma \in B_{n}$, which we call its signed index permutation, 
that satisfies the following properties:
\begin{enumerate}
\item $q_{|\sigma(1)|}\ge q_{|\sigma(2)|}\ge \cdots\ge q_{|\sigma(n)|}$, 
\item if $0<i<j$ and $q_{|\sigma(i)|}=q_{|\sigma(j)|}$ then
$\sigma(i)<\sigma(i+1)<\cdots<\sigma(j)$, and 
\item $q_{|\sigma(i)|}$ is even if and only if $\sigma(i)>0$. 
\end{enumerate}
In other words, the signed permutation $\sigma$ is the unique 
element in $B_{n}$ whose signs are determined by the parity of the 
$q_{i}$'s and that orders the elements in the sequence 
$q=(q_{1},\dots,q_{n})$ in decreasing way breaking ties 
from left to right for even values of $q_{i}$ and from right to left 
for odd values of $q_{i}$.

\medskip

\noindent{\bf{Example:}} Suppose
 $m(\x,\y)=x^{7}_{1}x^{6}_{2}x^{6}_{3}x^{5}_{4}x^{5}_{5}x^{3}_{6}
y^{3}_{1}y^{8}_{2}y^{6}_{3}y^{3}_{4}y^{5}_{5}y_{6}^{5}$. Then the 
signed index permutation associated to $m(\x,\y)$ is 
$\sigma=[2,3,-6,-5,-4,-1]$.
 
\subsection{Exponent decomposition}
Suppose that $m(\x,\y)=x^{p}y^{q}$ is a monomial in $\O_{n}$ 
and let $\sigma\in B_{n}$ be the signed index permutation 
associated to $m(\x,\y)$ as explained above. We can use the 
signed permutation $\sigma$ to obtain a decomposition of the 
sequences $p$ and $q$ as we explain next. We start by 
decomposing $q$.  

\begin{property}\label{Property q}
The sequence  $\{q_{|\sigma(i)|}-f_{i}(\sigma)\}_{i=1}^{n}$  
is a decreasing sequence of non-negative even integers.
\end{property} 
\Proof
By definition and property (3) of the signed index permutation 
$\sigma$ we have
\begin{align*}
q_{|\sigma(i)|}-f_{i}(\sigma)
&=q_{|\sigma(i)|}-2d_{i}(\sigma)-\varepsilon_{i}(\sigma)\\
&\equiv q_{|\sigma(i)|}-\varepsilon_{i}(\sigma) \  (\text{mod } 2)\\
&\equiv 0 \  (\text{mod } 2).
\end{align*}
This proves that $q_{|\sigma(i)|}-f_{i}(\sigma)$ is even for all 
$1\le i\le n$. Note that $q_{|\sigma(n)|}-f_{n}(\sigma)
=q_{|\sigma(n)|}-\varepsilon_{n}(\sigma)\ge 0$ since 
$q_{|\sigma(n)|}\ge 0$ and $q_{|\sigma(n)|}$ and 
$\varepsilon_{i}(\sigma)$  have the same parity. It remains 
to prove $q_{|\sigma(i)|}-f_{i}(\sigma)\ge 
q_{|\sigma(i+1)|}-f_{i+1}(\sigma)$ for all $1\le i\le n-1$.  
For this  we consider the following cases.

$\bullet$ Case 1. Suppose that $\sigma(i)<\sigma(i+1)$.  
This implies $d_{i}(\sigma)=d_{i+1}(\sigma)$. Thus in this 
case we need to prove that 
$q_{|\sigma(i)|}-\varepsilon_{i}(\sigma)\ge 
q_{|\sigma(i+1)|}-\varepsilon_{i+1}(\sigma)$. Since 
$q_{|\sigma(i)|}\ge q_{|\sigma(i+1)|}$ 
the only case we need to inspect is the case 
$\varepsilon_{i}(\sigma)=1$ and $\varepsilon_{i+1}(\sigma)=0$. 
However, under this assumption $q_{|\sigma(i)|}$ is odd and 
$q_{|\sigma(i+1)|}$ is even 
and thus $q_{|\sigma(i)|}-1\ge q_{|\sigma(i+1)|}$.
 
$\bullet$ Case 2.  Suppose that $\sigma(i)>\sigma(i+1)$ 
and $\varepsilon_{i}(\sigma)\ne \varepsilon_{i+1}(\sigma)$.
This implies  $\varepsilon_{i}(\sigma)=0$  and 
$\varepsilon_{i+1}(\sigma)=1$. We have 
$d_{i}(\sigma)=d_{i+1}(\sigma)+1$. In this case we need 
to show that $q_{|\sigma(i)|}-1\ge q_{|\sigma(i+1)|}$. 
Note that $q_{|\sigma(i)|}$  must be even and $q_{|\sigma(i+1)|}$ 
must be odd and by property (1) of the signed index permutation 
$q_{|\sigma(i)|}\ge q_{|\sigma(i+1)|}$. Therefore 
$q_{|\sigma(i)|}> q_{|\sigma(i+1)|}$ which means
$q_{|\sigma(i)|}-1\ge q_{|\sigma(i+1)|}$

$\bullet$ Case 3.  Suppose that $\sigma(i)>\sigma(i+1)$ 
and $\varepsilon_{i}(\sigma)=\varepsilon_{i+1}(\sigma)$.
Then  $d_{i}(\sigma)=d_{i+1}(\sigma)+1$. In this case we need 
to show that $q_{|\sigma(i)|}-2\ge q_{|\sigma(i+1)|}$. Since 
$\varepsilon_{i}(\sigma)=\varepsilon_{i+1}(\sigma)$, 
then $q_{|\sigma(i)|}$ and $q_{|\sigma(i+1)|}$ must have the 
same parity and by condition (1) of the 
signed index permutation we have 
$q_{|\sigma(i)|}\ge q_{|\sigma(i+1)|}$. Thus we only need to 
prove that $q_{|\sigma(i)|}> q_{|\sigma(i+1)|}$. Assume by 
contradiction that $q_{|\sigma(i)|}=q_{|\sigma(i+1)|}$. 
Using condition (2) of the signed index permutation 
we conclude $\sigma(i)<\sigma(i+1)$ which contradicts 
our original assumption.
\qed

\medskip

By the previous property, for every $1\le i\le n$ 
we can find a non-negative number 
$\mu_{|\sigma(i)|}$  such that 
$q_{|\sigma(i)|}=2\mu_{|\sigma(i)|}+f_{i}(\sigma)$. Define 
$\gamma_{|\sigma(i)|}:=f_{i}(\sigma)$ so that 
$\gamma_{i}=f_{|\sigma^{-1}(i)|}(\sigma)$. 

\begin{proposition}\label{decom q}
The sequences $\gamma=(\gamma_{1},\dots, \gamma_{n})$ 
and $\mu=(\mu_{1},\dots,\mu_{n})$ are sequences of non-negative 
integers that satisfy the following 
properties:
\begin{enumerate}
\item $q=2\mu+\gamma$, 
\item $\mu_{|\sigma(1)|}\ge\mu_{|\sigma(2)|}
\ge\cdots\ge \mu_{|\sigma(n)|}$,
\item $\gamma_{|\sigma(1)|}\ge\gamma_{|\sigma(2)|}
\ge\cdots\ge \gamma_{|\sigma(n)|}$,
\item if $0<i<j\le n$ and $\gamma_{i}=\gamma_{j}$ then 
$\s(q_{i})\ge \s(q_{j})$.
\end{enumerate}
\end{proposition}
\Proof
Properties (1)--(3) follow directly from the construction and 
from the properties of the numbers $f_{i}(\sigma)$. Suppose  
now that $0<i<j\le n$ and $\gamma_{i}=\gamma_{j}$. This means 
$f_{|\sigma^{-1}(i)|}(\sigma)=f_{|\sigma^{-1}(j)|}(\sigma)$;  that is, 
$2d_{|\sigma^{-1}(i)|}(\sigma)+\varepsilon_{|\sigma^{-1}(i)|}(\sigma)=
2d_{|\sigma^{-1}(j)|}(\sigma)+\varepsilon_{|\sigma^{-1}(j)|}(\sigma)$. 
It follows that $\varepsilon_{|\sigma^{-1}(i)|}(\sigma)
=\varepsilon_{|\sigma^{-1}(j)|}(\sigma)$ and 
$d_{|\sigma^{-1}(i)|}(\sigma)=d_{|\sigma^{-1}(j)|}(\sigma)$. By definition 
\begin{align*}
q_{i}&=\gamma_{i}+2\mu_{i}=2(\mu_{i}+2d_{|\sigma^{-1}(i)|}(\sigma))
+\varepsilon_{|\sigma^{-1}(i)|}(\sigma),\\
q_{j}&=\gamma_{j}+2\mu_{j}=2(\mu_{j}+2d_{|\sigma^{-1}(j)|}(\sigma))
+\varepsilon_{|\sigma^{-1}(j)|}(\sigma).
\end{align*}
In particular we conclude  that $q_{i}$ and $q_{j}$ have the same parity. 
We need to consider two cases according to the parity of these numbers. 
Suppose first that $q_{i}$ and $q_{j}$ 
are even. Let $k=|\sigma^{-1}(i)|$ and $l=|\sigma^{-1}(j)|$. Since 
$q_{i}$ and $q_{j}$ are even, then 
$\varepsilon_{k}(\sigma)=\varepsilon_{l}(\sigma)=0$ and this implies 
that $\sigma^{-1}(i),\sigma^{-1}(j)>0$; that is ,
$k=\sigma^{-1}(i)$ and $l=\sigma^{-1}(j)$. Let's show that 
$k<l$. Assume by contradiction that $l<k$. Since 
$d_{l}(\sigma)=d_{k}(\sigma)$ and $l < k$, then by property (2) of 
the numbers $d_{i}(\sigma)$ it follows
\[
j=\sigma(l)<\sigma(l+1)<\cdots<\sigma(k)=i
\]
which contradicts the assumption  $i<j$. Therefore
$0<\sigma^{-1}(i)<\sigma^{-1}(j)$ and  by property 
(1) of the signed index permutation we conclude 
$q_{i}=q_{|\sigma(\sigma^{-1}(i))|}\ge 
q_{|\sigma(\sigma^{-1}(j))|}=q_{j}$. The case where 
$q_{i}$ and $q_{j}$ are odd is handled in a similar way.
 \qed
 
 \medskip

Next we obtain a similar decomposition for the sequence 
$p=(p_{1},\dots, p_{n})$. To start assume 
$m(\x,\y)=x^{p}y^{q}\in \O_{n}$.  This implies that
 $p_{i}+q_{i}$ is even for all $1\le i\le n$. On the other hand,  
$\varepsilon_{i}(\sigma^{-1})=0$ if and only if 
$k:=\sigma^{-1}(i)>0$ and this is the case if and only if 
$q_{i}=q_{\sigma(k)}$ is even. We conclude 
that if $m(\x,\y)\in \O_{n}$ then for  every $1\le i\le n$
\[
p_{i}\equiv q_{i}\equiv \varepsilon_{i}(\sigma^{-1}) \ (\text{mod } 2).
\]

\begin{property}\label{Property p}
Suppose that $m(\x,\y)=x^{p}y^{q}\in \O_{n}$. Then 
$\{p_{i}-f_{i}(\sigma^{-1})\}_{i=1}^{n}$  
is a decreasing sequence of non-negative even integers.
\end{property}
\Proof
By the above comment 
\begin{align*}
p_{i}-f_{i}(\sigma^{-1})
&=p_{i}-2d_{i}(\sigma^{-1})-\varepsilon_{i}(\sigma^{-1})\\
&\equiv q_{i}-\varepsilon_{i}(\sigma^{-1}) \  (\text{mod } 2)\\
&\equiv 0 \  (\text{mod } 2).
\end{align*}
This proves that $p_{i}-f_{i}(\sigma^{-1})$ is even for all $1\le i\le n$. 
On the other hand, note that $p_{n}-f_{n}(\sigma^{-1})
=p_{n}-\varepsilon_{n}(\sigma^{-1})\ge 0$ 
because $p_{n}$ is odd if and only if $\varepsilon_{n}(\sigma^{-1})=1$. 
We are left to prove that 
$p_{i}-f_{i}(\sigma^{-1})\ge p_{i+1}-f_{i+1}(\sigma^{-1})$ for all 
$1\le i\le n-1$.  For this  we consider 
the following cases.

$\bullet$ Case 1. Suppose that $\sigma^{-1}(i)<\sigma^{-1}(i+1)$.  
This implies $d_{i}(\sigma^{-1})=d_{i+1}(\sigma^{-1})$. Thus in this 
case we need to prove that 
$p_{i}-\varepsilon_{i}(\sigma^{-1})\ge 
p_{i+1}-\varepsilon_{i+1}(\sigma^{-1})$. Since $p_{i}\ge p_{i+1}$ 
the only case we need to inspect is the case 
$\varepsilon_{i}(\sigma^{-1})=1$ and 
$\varepsilon_{i+1}(\sigma^{-1})=0$. 
However, under this assumption $p_{i}$ is odd and $p_{i+1}$ is even 
and thus $p_{i}-1\ge p_{i+1}$ as $p_{i}\ge p_{i+1}$.
 
$\bullet$ Case 2.  Suppose that $\sigma^{-1}(i)>\sigma^{-1}(i+1)$ 
and $\varepsilon_{i}(\sigma^{-1})\ne \varepsilon_{i+1}(\sigma^{-1})$.
This is only possible if $\varepsilon_{i}(\sigma^{-1})=0$  and 
$\varepsilon_{i+1}(\sigma^{-1})=1$. Then 
$i\in Des(\sigma^{-1})$ and 
$d_{i}(\sigma^{-1})=d_{i+1}(\sigma^{-1})+1$. In this case we need 
to show that $p_{i}-1\ge p_{i+1}$. 
Note that $p_{i}$  must be even and $p_{i+1}$ must be odd 
and by assumption $p_{i}\ge p_{i+1}$. Therefore $p_{i}>p_{i+1}$ 
and thus $p_{i}-1\ge p_{i+1}$ as desired.
 
$\bullet$ Case 3.   Suppose that $\sigma^{-1}(i)>\sigma^{-1}(i+1)$ 
and $\varepsilon_{i}(\sigma^{-1})=\varepsilon_{i+1}(\sigma^{-1})$. 
Then $i\in Des(\sigma^{-1})$ and therefore 
$d_{i}(\sigma^{-1})=d_{i+1}(\sigma^{-1})+1$. In this case we need 
to show that $p_{i}-2\ge p_{i+1}$. Note that 
$\varepsilon_{i}(\sigma^{-1})=\varepsilon_{i+1}(\sigma^{-1})$ implies 
that $p_{i}$ and $p_{i+1}$ have the same parity. We know that 
$p_{i}\ge p_{i+1}$. Therefore we only need to prove that 
$p_{i}>p_{i+1}$. Assume by contradiction that $p_{i}=p_{i+1}$. 
By assumption $m(\x,\y)=x^{p}y^{q}\in \O_{n}$ and $p_{i}=p_{i+1}$. 
This implies  $\s(q_{i})\ge\s(q_{i+1})$. 
Let $k=|\sigma^{-1}(i)|$ and $l=|\sigma^{-1}(i+1)|$. If  
$\varepsilon_{i}(\sigma^{-1})=\varepsilon_{i+1}(\sigma^{-1})=0$ then 
we obtain $0<l<k$ and therefore we conclude 
$q_{i+1}=q_{|\sigma(l)|}\ge q_{|\sigma(k)|}=q_{i}$ 
by property  (1) of the signed index permutation $\sigma$. Note that 
$q_{i}$ and $q_{i+1}$ must be even since they have the same parity 
as  $\varepsilon_{i}(\sigma^{-1})=\varepsilon_{i+1}(\sigma^{-1})=0$. 
Thus $q_{i}=\s(q_{i})\ge\s(q_{i+1})=q_{i+1}$ and in turn
$q_{i+1}=q_{|\sigma(l)|}= q_{|\sigma(k)|}=q_{i}$. This together with 
property (2) of the signed index permutation $\sigma$ imply 
$i+1=\sigma(l)<\sigma(k)=i$ which is a contradiction. Suppose now  
that $\varepsilon_{i}(\sigma^{-1})=\varepsilon_{i+1}(\sigma^{-1})=1$.  
Then $k=-\sigma^{-1}(i)$ and $l=-\sigma^{-1}(i+1)$ and by 
assumption $\sigma^{-1}(i)>\sigma^{-1}(i+1)$. Thus $0>-k>-l$; 
that is, $0<k<l$. Using property (1) of the signed index permutation 
we conclude $q_{i}=q_{|\sigma(k)|}\ge q_{|\sigma(l)|}=q_{i+1}$. 
Under the given assumptions $q_{i}$ and $q_{i+1}$ must be odd 
and $\s(q_{i})\ge\s(q_{i+1})$; that is, $q_{i}\le q_{i+1}$. Again 
we conclude that $q_{i}=q_{|\sigma(k)|}= q_{|\sigma(l)|}=q_{i+1}$.
Since $0<k<l$ using property (2)  as before  
we conclude that  $-i=\sigma(k)<\sigma(l)=-i-1$ deriving a 
contradiction in either case.
 \qed

\medskip

By the previous property, for every $1\le i\le n$ we can find 
a non-negative integer $\nu_{i}$ such that 
$2\nu_{i}=p_{i}-f_{i}(\sigma^{-1})$; that is, 
$p_{i}=2\nu_{i}+f_{i}(\sigma^{-1})$. Define 
$\delta_{i}:=f_{i}(\sigma^{-1})$.

\begin{proposition}\label{decom p} 
The sequences $\delta=(\delta_{1},\dots, \delta_{n})$ 
and $\nu=(\nu_{1},\dots,\nu_{n})$ are sequences of non-negative 
integers that satisfy the following 
properties:
\begin{enumerate}
\item $p=2\nu+\delta$, 
\item $\nu_{1}\ge\nu_{2}
\ge\cdots\ge \nu_{n}$,
\item $\delta_{1}\ge\delta_{2}
\ge\cdots\ge \delta_{n}$,
\item if $0<i<j\le n$ and $\delta_{i}=\delta_{j}$ then 
$\s(q_{i})\ge \s(q_{j})$.
\end{enumerate}
\end{proposition}
\Proof 
Properties (1)--(3) follow directly from the construction, 
Property \ref{Property p} and  the properties of the numbers 
$f_{i}(\sigma^{-1})$ . Suppose  
now that $0<i<j\le n$ and $\delta_{i}=\delta_{j}$. This means 
$f_{i}(\sigma^{-1})=f_{j}(\sigma^{-1})$ and in turn 
$d_{i}(\sigma^{-1})=d_{j}(\sigma^{-1})$ and 
$\varepsilon_{i}(\sigma^{-1})=\varepsilon_{j}(\sigma^{-1})$.  
Since $i<j$ by property (2) of the numbers $d_{i}(\sigma^{-1})$ 
it follows that $\sigma^{-1}(i)<\sigma^{-1}(j)$.  Since 
$\varepsilon_{i}(\sigma^{-1})=\varepsilon_{j}(\sigma^{-1})$ then 
$q_{i}$ and $q_{j}$ have the same parity. Assume $q_{i}$ 
and $q_{j}$ are both even, then 
$0<\sigma^{-1}(i)<\sigma^{-1}(j)$ and by property 
(1) of the signed index permutation we have 
\[
q_{i}=q_{|\sigma(\sigma^{-1}(i))|}\ge 
q_{|\sigma(\sigma^{-1}(j))|}=q_{j}.
\]
Similarly, if $q_{i}$ and $q_{j}$ are odd, then 
$\sigma^{-1}(i)<\sigma^{-1}(j)<0$ and thus 
$-\sigma^{-1}(i)>-\sigma^{-1}(j)>0$. Again by property 
(1) of the number signed index permutation we have 
\[
q_{i}=q_{|\sigma(\sigma^{-1}(i))|}\le 
q_{|\sigma(\sigma^{-1}(j))|}=q_{j}.
\]
In either case we obtain $\s(q_{i})\ge \s(q_{j})$.
\qed

\subsection{Monomial decomposition}  
Suppose that $m(\x,\y)=x^{p}y^{q}\in \O_{n}$. We can 
decompose $p=2\nu+\delta$ and $q=2\mu+\gamma$ using 
Propositions \ref{decom p} and \ref{decom q}. Let $\sigma\in B_{n}$  
be the index permutation associated to the monomial $m(\x,\y)$. 
By definition $\delta_{i}=f_{i}(\sigma^{-1})$ and
$\gamma_{i}=f_{|\sigma^{-1}(i)|}(\sigma)$. Therefore
\[
c_{\sigma}=\prod_{i=1}^{n}x_{i}^{f_{i}(\sigma^{-1})}
y_{i}^{f_{|\sigma^{-1}(i)|}(\sigma)}=x^{\delta}y^{\gamma}.
\] 
This means that we have a decomposition
\[
m(\x,\y)=x^{p}y^{q}=x^{2\nu+\delta}y^{2\mu+\gamma}
=x^{2\nu}y^{2\mu}c_{\sigma}.
\]
Given a sequence of integers $i=(i_{1},\dots, i_{n})$ 
let $\Sigma_{n}(i)$  denote the stabilizer of $i$ under 
the  permutation action of $\Sigma_{n}$; that is, 
$\Sigma_{n}(i)$ is the subgroup of elements in  $\Sigma_{n}$
that fix $i$. Define
\[
m_{2\nu}(\x)=
\sum_{[\alpha]\in\Sigma_{n}/\Sigma_{n}(\nu)}x^{2\alpha(\nu)}, 
\ \text{ and } \
m_{2\mu}(\y)=
\sum_{[\beta]\in\Sigma_{n}/\Sigma_{n}(\mu)}y^{2\beta(\mu)}.
\]
By definition $m_{2\nu}(\x)$ and 
$m_{2\mu}(\y)$ are symmetric polynomials on the variables 
$x_{1}^{2},\dots, x_{n}^{2}$ and $y_{1}^{2},\dots, y_{n}^{2}$ 
respectively. In particular 
$m_{2\nu}(\x)m_{2\mu}(\y)\in R^{B_{n}}$.  
In the same way as in \cite[Proposition 3.2]{Allen} the polynomial 
$m_{2\nu}(\x)m_{2\mu}(\y)\rho(c_{\sigma})$ can be decomposed 
as we prove next.

\begin{theorem}
Suppose that $m(\x,\y)\in \O_{n}$ and let $\sigma\in B_{n}$ 
be the corresponding singed index permutation. Then 
\[
m_{2\nu}(\x)m_{2\mu}(\y)\rho(c_{\sigma})= k_{m(\x,y)}\rho(m(\x,\y))
+\sum_{n(\x,\y)\precneqq m(\x,\y)}k_{n(\x,\y)}\rho(n(\x,\y)).
\]
Here $n(\x,y)\in \O_{n}$ 
runs through the collection of ordered monomials 
with same total degree as $m(\x,y)$ with 
$n(\x,\y)\precneqq m(\x,\y)$, and
$k_{m(\x,y)}>0$, $k_{n(\x,\y)}$ are constants. 
\end{theorem}
\Proof 
Using the definitions we have
\begin{align*}
m_{2\nu}(\x)m_{2\mu}(\y)\rho(c_{\sigma})
&=\rho(m_{2\nu}(\x)m_{2\mu}(\y)c_{\sigma})=
\sum_{[\alpha]\in\Sigma_{n}/\Sigma_{n}(\nu)}
\sum_{ [\beta]\in\Sigma_{n}/\Sigma_{n}(\mu)}
\rho(x^{2\alpha\nu}y^{2\beta\mu}c_{\sigma})\\
&=\sum_{[\alpha]\in\Sigma_{n}/\Sigma_{n}(\nu)}
\sum_{ [\beta]\in\Sigma_{n}/\Sigma_{n}(\mu)}
\rho(x^{2\alpha\nu+\delta}y^{2\beta\mu+\gamma}).
\end{align*}
Fix $\alpha, \beta\in \Sigma$  and consider 
$n(\x,y)\in \O_{n}$ the unique monomial such that 
$\rho(n(\x,\y))=\rho(x^{2\alpha\nu+\delta}
y^{2\beta\mu+\gamma})$. 
To prove the theorem we need to prove that 
$n(\x,\y)\preccurlyeq m(\x,\y)$.  
Let $[\Sigma_{n}(\delta)]$ denote the image of 
$\Sigma_{n}(\delta)$ in $\Sigma_{n}/\Sigma_{n}(\nu)$ under 
the natural map. By parts (2) and (3)  of Proposition \ref{decom p} 
it follows that if 
$[\alpha]\in [\Sigma_{n}(\delta)]$ then 
$\o(2\alpha\nu+\delta)=\o(p)$, and if 
$[\alpha] \notin [\Sigma_{n}(\delta)]$ 
then $\o(2\alpha\nu+\delta)<_{\ell}\o(p)$.  Similarly, 
if $[\Sigma_{n}(\gamma)]$ denote the image of 
$\Sigma_{n}(\gamma)$ in  $\Sigma_{n}/\Sigma_{n}(\mu)$
then by parts (2) and (3)  of Proposition \ref{decom q} 
it follows that if 
$[\beta]\in [\Sigma_{n}(\gamma)]$ then 
$\o(2\beta\mu+\gamma)=\o(q)$, and if 
$[\beta] \notin [\Sigma_{n}(\gamma)]$ 
then $\o(2\beta\mu+\gamma)<_{\ell}\o(q)$.  
With this in mind we have the following cases.

$\bullet$ Case 1.  Suppose that $[\alpha]\notin [\Sigma_{n}(\delta)] $ 
or $[\beta]\notin [\Sigma_{n}(\gamma)] $.
If $[\alpha]\notin [\Sigma_{n}(\delta)]$ then by the previous comment
$\o(2\alpha\nu+\delta)<_{\ell}\o(p)$ and if 
$[\alpha]\in [\Sigma_{n}(\delta)] $ but  
$\beta\notin [\Sigma_{n}(\gamma)] $ 
then $\o(2\alpha\nu+\delta)=\o(p)$ but 
$\o(2\beta\mu+\gamma)<_{\ell}\o(q)$. In either case we 
conclude
\[
\o(n(\x,\y))=(\o(2\alpha\nu+\delta),\o(2\beta\mu+\gamma)) 
<_{\ell}(\o(p),\o(q))
=\o(m(\x,\y)).
\]
It follows that in this case $n(\x,\y)\precneqq m(\x,\y)$.
 
$\bullet$ Case 2. Suppose that  $[\alpha]\in [\Sigma_{n}(\delta)] $ 
and $[\beta]\in [\Sigma_{n}(\gamma)] $. We can assume without loss 
of generality that 
$\alpha\in \Sigma_{n}(\delta) $ 
and $\beta\in \Sigma_{n}(\gamma)$. Then 
\[
x^{2\alpha\nu+\delta}y^{2\beta\mu+\gamma}
=x^{\alpha(2\nu+\delta)}y^{\beta(2\mu+\gamma)} 
=x^{\alpha(p)}y^{\beta(q)}.
\]
Therefore $n(\x,\y)=x^{p}y^{\pi \alpha^{-1}\beta(q)}$ for some 
$\pi \in \Sigma_{n}$. Note that the permutation 
$\pi$ has to stabilize $p$ and $\delta$, thus in particular 
$\pi\in \Sigma_{n}(\delta)$. In this case 
\[
\o(n(\x,\y))=\o(x^{\alpha(p)}y^{\beta(q)})=\o(m(\x,\y)).
\]
To prove that $n(\x.y)\preccurlyeq m(\x,\y)$ we need to show that 
$(p,\s(\pi \alpha^{-1}\beta(q)))\le_{\ell}(p,\s(q))$; that is, we need 
to prove that $\s(\pi \alpha^{-1}\beta(q))\le_{\ell} \s(q)$.  If 
$\s(\pi \alpha^{-1}\beta(q))=\s(q)$ there is nothing to prove. 
Assume that $\s(\pi \alpha^{-1}\beta(q))\ne \s(q)$ 
and let $1\le k\le n$ be the smallest integer such that 
$\s(q_{\pi \alpha^{-1}\beta(k)})\ne \s(q_{k})$.  We need to show 
that $\s(q_{\pi \alpha^{-1}\beta(k)})< \s(q_{k})$.  
Since $\beta\in \Sigma_{n}(\gamma)$, then by 
Proposition \ref{decom q} part (4)  we have that if 
$i<\beta(i)$ then $\s(q_{i})\ge \s(q_{\beta(i)})$. Similarly, 
since $\alpha,\pi\in \Sigma_{n}(\delta)$ by
Proposition \ref{decom p} part (4) whenever $i<\pi\alpha^{-1}(i)$ 
then $\s(q_{i})\ge \s(q_{\pi\alpha^{-1}(i)})$.  Using this 
we can see that $\s(q_{k})> \s(q_{\pi \alpha^{-1}\beta(k)})$. 

In either case we conclude that   $n(\x,\y)\preccurlyeq m(\x,\y)$.
\qed

\subsection{Main theorem} Finally we are ready to prove the main 
theorem of this article.

\begin{theorem}
Suppose that $n\ge 1$. Then the collection 
$\{\rho(c_{\sigma})\}_{\sigma\in B_{n}}$ forms a free basis of $M^{B_{n}}$ 
as a module over $R^{B_{n}}$.
\end{theorem}
\Proof
Let's show first that $\{\rho(c_{\sigma})\}_{\sigma\in B_{n}}$ 
generates $M^{B_{n}}$ as a module over $R^{B_{n}}$. It suffices to 
show that for every  $m(\x,\y)\in \O_{n}$ the polynomial 
$\rho(m(\x,\y))$ is generated by the different $\rho(c_{\sigma})$. 
Fix $m(\x,\y)\in \O_{n}$ and let $\sigma$ be the corresponding 
singed index permutation. By the previous theorem we have 
\[
m(\x,\y)=l_{m(\x,\y)}m_{2\nu}(\x)m_{2\mu}(\y)\rho(c_{\sigma})
+\sum_{n(\x,\y)\precneqq m(\x,\y)}l_{n(\x,\y)}n(\x,\y),
\]
for some constants $l_{n(\x,\y)}$ and monomials 
$n(\x,\y)\precneqq m(\x,\y)$ of same total degree. 
Iterating this process on the
monomials $n(\x,\y)\precneqq m(\x,\y)$ as many times as 
necessary 
we see that we can write  $m(\x,\y)$ as a 
linear combination of the elements 
$\{\rho(c_{\sigma})\}_{\sigma\in B_{n}}$ with coefficients 
in $R^{B_{n}}$. (Note that this 
process must terminate after finitely many stages as there 
are only finitely many monomials $n(\x,\y)$ of same 
total degree as $m(\x,\y)$). This proves 
that $\{\rho(c_{\sigma})\}_{\sigma\in B_{n}}$ generates 
$M^{B_{n}}$ as a $R^{B_{n}}$-module. On the other hand, 
note that $M^{B_{n}}$ 
is a bigraded ring over $\Q$ with 
\[
{\rm{bideg}}(x^{p}y^{q})=(|p|,|q|),
\]
where $|p|=p_{1}+\cdots+p_{n}$. With this 
graduation for every $\sigma\in B_{n}$ we 
have that  the polynomial 
$\rho(c_{\sigma})$ is homogeneous and 
\[
{\rm{bideg}}(\rho(c_{\sigma}))=(\fmaj(\sigma^{-1}),\fmaj(\sigma)).
\]
Let $P_{M^{B_{n}}}(s,t)$ denote 
the bigraded Hilbert series of the bigraded ring $M^{B_{n}}$. 
Using \cite[Theorem 3]{AR} we conclude 
that the series $P_{M^{B_{n}}}(s,t)$ 
is given by
\[
P_{M^{B_{n}}}(s,t)=\frac{\left(\sum_{\sigma\in B_{n}}
s^{\fmaj(\sigma^{-1})}t^{\fmaj(\sigma)}\right)}
{\prod_{i=1}^{n}(1-s^{2i})(1-t^{2i})}.
\]
This together with \cite[Theorem 1.4]{Allen} show that 
$\{\rho(c_{\sigma})\}_{\sigma\in B_{n}}$ is a free basis of 
$M^{B_{n}}$ as module over  $R^{B_{n}}$.
\qed

\medskip

As in the previous cases this theorem also has a 
geometric application. Let $B(2,Sp(n))$ be the geometric 
realization of the simplicial space 
obtained by considering the space of 
commuting $k$-tuples in $Sp(n)$, 
$B_{k}(2,Sp(n))=\Hom(\Z^{k},Sp(n))$. 
In \cite{AG} it is proved that the signed diagonal 
descent monomials can be used to obtain an explicit basis of 
$H^{*}(B(2,Sp(n));\Q)$ seen as a module over $H^{*}(BSp(n);\Q)$.

\end{document}